\documentclass[11pt,fleqn]{article}
\RequirePackage[OT1]{fontenc}
\RequirePackage{amsthm,amsmath}
\RequirePackage[round]{natbib}
\RequirePackage[colorlinks,citecolor=blue,urlcolor=blue]{hyperref}

\usepackage{amsmath,amsfonts,amssymb,euscript,mathrsfs}
\usepackage{slashbox}
\usepackage{float}
\usepackage[final]{pdfpages}

\usepackage{txfonts}
\usepackage[T1]{fontenc} 
\usepackage{fullpage}
\usepackage{graphicx}
\usepackage{dsfont}
\makeatletter
\def\captionof#1#2{{\def\@captype{#1}#2}}
\makeatother
\def\1{\mbox{\bf 1}}
\def\R{\mathbb{R}}

\def\N{\mathbb{N}}
\def\P{\mathbb{P}}
\def\E{\mathbb{E}}
\def\L{\mathbb{L}}

\def\R{\mathbb{R}}

\def\Z{\mathbb{Z}}

\def\H{\mathbb{H}}
\def\v{\mbox{Var\,}}

\def\c{\mbox{Cov}}

\newtheorem{theo}{Theorem}
\newtheorem{lem}{Lemma}
\newtheorem{prop}{Proposition}

\newtheorem{cor}{Corollary}
\newtheorem{Def/Prop}{Definition-Proposition}

\newcounter{exos}
\renewcommand\theexos{\arabic{exos}}

\newcounter{prob}
\renewcommand\theprob{\arabic{prob}}

\begin{document}
\author{Lionel Truquet \footnote{UMR 6625 CNRS IRMAR, University of Rennes 1, Campus de Beaulieu, F-35042 Rennes Cedex, France and}
\footnote{ENSAI, Campus de Ker-Lann, rue Blaise Pascal, BP 37203, 35172 Bruz cedex, France. {\it Email: lionel.truquet@ensai.fr}.}
 }
\title{Efficient semiparametric estimation in time-varying  regression models}
\date{}
\maketitle

\begin{abstract}
We study semiparametric inference in some linear regression models with time-varying coefficients, dependent regressors and dependent errors. This problem, which has been considered recently by \citet{ZW} under the functional dependence measure, is interesting for parsimony reasons or for testing stability of some coefficients in a linear regression model. 
In this paper, we propose a different procedure for estimating non time-varying parameters at the rate $\sqrt{n}$, in the spirit of the method introduced by \citet{Rob} for partially linear models. When the errors in the model are martingale differences, this approach can lead to more efficient estimates than the method considered in \citet{ZW}. For a time-varying AR process with exogenous covariates and conditionally Gaussian errors, we derive a notion of efficient information matrix from a convolution theorem adapted to triangular arrays. For  independent but non identically distributed Gaussian errors, we construct an asymptotically efficient estimator in a semiparametric sense.
\end{abstract}

\section{Introduction}
Regression models with time-varying parameters have been widely studied in the literature. This type of models includes standard regression models, some linear time series models such as the locally stationary AR process studied by \citet{Dahlhaus} or some time series models with covariates considered for instance in \citet{Chen}. 
In \citet{ZW}, the authors study a general linear regression model with time-varying coefficients and dependent errors. This model writes 
\begin{equation}\label{model1}
y_i=x_i'\beta(i/n)+e_i,\quad 1\leq i\leq n,\quad \E\left(e_i\vert x_i\right)=0,
\end{equation} 
where $\left(x_i\right)_{1\leq i\leq n}$ and $\left(e_i\right)_{1\leq i\leq n}$ are triangular arrays of dependent random vectors and 
$\beta$ is a smooth function which defines the time-varying coefficients in the regression model. 
Model (\ref{model1}) is quite general under the assumptions used in this paper. For instance, the regressors and the noise component can be stationary processes or locally stationary processes in the sense of \citet{Dahlhaus}.   
An interesting and natural question for model (\ref{model1}) is to reduce model complexity by assuming that some covariates are associated to a constant parameter. In \citet{ZW}, a testing procedure is given for detecting some stable (or non time-varying) coefficients. Under the null hypothesis, one faces to a semiparametric model of type
\begin{equation}\label{model2}
y_i=x_{1,i}'\beta_1+x_{2,i}'\beta_2(i/n)+e_i,\quad 1\leq i\leq n,$$
\end{equation} 
with $x_i=\left(x_{i,1}',x_{i,2}'\right)'$. 
Semiparametric models are often an interesting alternative when the number of regressors is important or when the sample size is moderate. Moreover, semiparametric models are often a first step to modify a parametric model by taking in account of a particular 
inhomogeneity in the data. For instance, \citet{Gao} investigated a regression model in which the time-varying trend is the single nonparametric component. In the same spirit, in \citet{Truquet2016}, we studied semiparametric versions of locally stationary ARCH models with a particular emphasis on the case of a time-varying unconditional variance and non time-varying lag parameters. The square of the ARCH process being a weak AR process, estimation of the ARCH coefficients is a particular case of statistical inference in model (\ref{model1}). 

As shown in \citet{ZW}, under some general assumptions, non time-varying parameters in model (\ref{model1}) can be estimated at the rate $\sqrt{n}$. The technique used in their paper consists in averaging the nonparametric estimator of the full vector $\beta(\cdot)$ 
for the coordinates corresponding to the stable parameters.  

In this paper, we generalize the approach studied in \citet{Truquet2016} for estimating stable parameters in ARCH models to model (\ref{model1}). Our method is based on the idea originally used by \citet{Rob} or \citet{SP} for partially linear models and consists in removing the nonparametric part of the regression function by using partial regressions.
This approach is also known as profile least squares method in the literature.
 Under some general assumptions, we show that $\beta_1$ can be estimated at the rate $\sqrt{n}$. But the main contribution of this paper addresses the question of efficiency. 
When the noise process $\left(e_i\right)_{1\leq i\leq n}$ is a martingale difference, we show that both methods requires the estimation of the conditional variance of the noise to improve efficiency. 
However, the lower bound for the asymptotic variance obtained by the averaging method of \cite{ZW} is larger than the method based on partial regressions. Moreover if the noise is conditionally Gaussian and the regressors exogenous (except possible lag values of $y_i$), the lower bound obtained with our approach can be interpreted as a lower bound in semiparametric estimation. Estimating the conditional variance and using a plug-in method to get an efficient estimate of $\beta_1$  is a very difficult task in general. This is why we only consider the case of independent errors with a time-varying unconditional variance to construct an asymptotically optimal estimator which is also efficient in the Gaussian case.

Let us mention that our conclusions for efficiency were already known for the case of a linear regression model with semi-varying coefficients, a model in which some of the coefficients depend on a covariate and the others are constant. See \citet{Fan1} for a survey of
regression models with varying coefficients and for a discussion of the non efficiency of the averaging method. The case of time-varying coefficients could be seen as a particular case of a semi-varying regression models with the time as a covariate. 
However, this particular case requires a special attention due to the dependence between the observations and the asymptotic with triangular arrays which lead to more technical arguments. Moreover, it is not possible to use some available results for the semiparametric efficiency such as that of \citet{Chamb}. Indeed, working with triangular arrays and the infill asymptotic inherent to the local stationarity assumptions requires a special attention.

In this paper, we will not study the nonparametric estimation of $\beta_2(\cdot)$ because using a plug-in of $\beta_1$ with a $\sqrt{n}-$consistent estimator in (\ref{model2}) is asymptotically equivalent to the nonparametric estimation of parameter $\beta_2$ with known coefficients $\beta_1$. The latter problem has been extensively studied in the literature. Then the most challenging problem is the estimation of the parametric part of the regression function, a problem we focus on in this paper.

The paper is organized as follows. In Section $2$, we introduce our model assumptions and we study the asymptotic properties of our estimators. In Section $3$, we consider  an autoregressive process with exogenous covariates and we give a discussion about asymptotic semiparametric efficiency when the errors are conditionally Gaussian.
Proofs of our two main results are postponed to the last section of the paper.

\section{Semiparametric estimation in time-varying regression models}
Let us introduce the model considered in this paper which is quite similar to the regression model studied in \citet{ZW}. We assume that
\begin{equation}\label{modele}
y_i=x_{1,i}'\beta_1+x_{2,i}'\beta_2(i/n)+e_i,\quad 1\leq i\leq n,
\end{equation}
where $e_i$ is a random variable taking values in $\R$, $x_{1,i}$ and $x_{2,i}$ are two random vectors of size $p_1$ and $p_2$ respectively and $\beta_1\in \R^{p_1}$, $\beta_2$ is a function defined on the unit interval $[0,1]$ and taking values in $\R^{p_2}$. In the sequel, we set $p=p_1+p_2$.
Now we introduce the dependence structure of $\left(x_{i,1}',x_{i,2}',e_i\right)'$ as well as the functional dependence measure used in \citet{ZW}. 
Let $\left(\xi_k\right)_{k\in \Z}$ a sequence of i.i.d random variables taking values in an arbitrary measurable space $\left(E,\mathcal{E}\right)$.
For $k\in\Z$, we set $\mathcal{F}_k=\sigma\left(\xi_j:j\leq k\right).$
 We also assume that there exist measurable functions 
$G_{n,i}:\left(\R^p\right)^{\N}\rightarrow \R^p$ and $H_{n,i}:\R^{\N}\rightarrow\R$ such that
$$x_i=\begin{pmatrix}x_{1,i}\\x_{2,i}\end{pmatrix}=G_{n,i}\left(\mathcal{F}_i\right),\quad e_i=H_{n,i}\left(\mathcal{F}_i\right).$$
If $\mathcal{F}_{0,k}=\sigma\left(\xi_0,\ldots,\xi_{-k+1},\xi'_{-k},\xi_{-k-1},\ldots\right)$ for $k\in\N$,
we set 
$$\delta_{k,q}(G)=\sup_{n\geq 1}\max_{1\leq i\leq n}\Vert G_{n,i}\left(\mathcal{F}_0\right)-G_{n,i}\left(\mathcal{F}_{0,k}\right)\Vert_q,$$
$$\Theta_{m,q}(G)=\sum_{k\geq m}\delta_{k,q}(G).$$
In the sequel, we will often use some conditions of type $\Theta_{0,q}(G)<\infty$. 
In this case, it will be assumed implicitly that $\max_{n\geq 1,1\leq i\leq n}\Vert G_{n,i}\Vert_q<\infty$.
Note that model (\ref{modele}) is quite general because it contains standard regression models, the autoregressive processes studied in \citet{ZW} and more generally autoregressive processes with  exogenous covariates. The latter case will be studied in Section $3$.  
The assumptions we will use in the sequel allow to consider locally stationary processes for the covariates and the noise component. This includes the moving averages considered in \citet{Dahlhaus}, the ARCH processes considered in \citet{DR} or other Markov processes studied in \citet{SR} and more recently in \citet{DaWu}. 
In particular the noise component $e_i$ can be conditionally heteroscedastic.

\subsection{Estimation of $\beta_1$}
Throughout this paper, $K$ will denote a probability density supported on $[-1,1]$ and Lipschitz continuous.
We set for $1\leq i,j\leq n$, 
$$k_{i,j}=\frac{\frac{1}{nb}K\left(\frac{i-j}{nb}\right)}{\frac{1}{nb}\sum_{\ell=1}^n K\left(\frac{i-\ell}{nb}\right)}.$$
We set
$$s_{1,i}=\E\left(x_{2,i}y_i\right),\quad s_{2,i}=\E\left(x_{2,i}x_{1,i}'\right),\quad s_{3,i}=\E\left(x_{2,i}x_{2,i}'\right),$$
$$q_{1,i}=s_{3,i}^{-1}s_{1,i},\quad q_{2,i}=s_{3,i}^{-1}s_{2,i},$$
$$\widetilde{y}_i=y_i-x_{2,i}'q_{1,i},\quad \widetilde{x}_{1,i}=x_{1,i}-q_{2,i}'x_{2,i}.$$
We have 
$$s_{1,i}=s_{2,i}\beta_1+s_{3,i}\beta_2(i/n).$$
This yields to the relations 
$$q_{1,i}=q_{2,i}\beta_1+\beta_2(i/n)$$
and 
$$\widetilde{y}_i=\widetilde{x}_{1,i}'\beta_1+e_i.$$
The quantities $\widetilde{y}_i$ and $\widetilde{x}_{1,i}$ will be estimated by 
$$\hat{y}_i=y_i-x_{2,i}'\hat{q}_{1,i},\quad \hat{x}_{1,i}=x_{1,i}-\hat{q}_{2,i}'x_{2,i},$$
where 
$$\hat{s}_{1,i}=\sum_{j=1}^n k_{i,j}x_{2,j}y_j,\quad \hat{s}_{2,i}=\sum_{j=1}^nk_{i,j}x_{2,j}x_{1,j}',\quad \hat{s}_{3,i}=\sum_{j=1}^n k_{i,j}x_{2,i}x_{2,i}',$$
$$\hat{q}_{1,i}=\hat{s}_{3,i}^{-1}\hat{s}_{1,i},\quad \hat{q}_{2,i}=\hat{s}_{3,i}^{-1}\hat{s}_{2,i}.$$
We now define the estimator of parameter $\beta_1$. 
$$\hat{\beta}_1=\left(\sum_{i=1}^n \hat{x}_{1,i}\hat{x}_{1,i}'\right)^{-1}\sum_{i=1}^n \hat{x}_{1,i}\hat{y}_i.$$
Our estimator can be interpreted in two ways. 
\begin{itemize}
\item
The first interpretation concerns partial regressions. The terms $x_{2,i}'q_{1,i}$ and $x_{2,i}'q_{2,i}$ can be interpreted as the orthogonal projection of $y_i$ and $x_{1,i}$ respectively onto the linear $\L^2$ space generated by the components of $x_{2,i}$. Then $\widetilde{y}_i$ and $\widetilde{x}_{1,i}$
are the residuals components of these two regressions and a new regression is made for estimating $\beta_1$. In practice, it is necessary to estimate non parametrically these quantities by $\hat{y}_i$ and $\hat{x}_{1,i}$ respectively. Then our estimator has the same interpretation than in the partially linear models studied in \citet{Rob} or \citet{SP}, except that the response $y_i$ and the regressors $x_{1,i}$ are projected onto the subspace of linear functions in $x_{2,i}$.

\item
The second interpretation is the profile least squares method studied for instance in \citet{Fan2} for semi-varying coefficients models. The quantity $\hat{q}_{1,i}-\hat{q}_{1,i}\beta_1$ is the Nadaraya-Watson estimator of the function $\beta_2$ for the model 
$$y_i-x_{1,i}'\beta_1=x_{2,i}'\beta_2(i/n)+e_i.$$
Then $\hat{\beta}_1$ is defined as a minimizer of
$$\alpha\mapsto \sum_{i=1}^n \left(y_i-x_{1,i}'\beta_1-x_{2,i}'\left[\hat{q}_{1,i}-\hat{q}_{1,i}\beta_1\right]\right)^2.$$
On can also use local polynomials in the first step if more regularity is assumed for $\beta_2$, see the remark after Theorem \ref{central}. 

\end{itemize}

Let $d$ be any positive integer. For a vector $v$ of $\R^d$ or a square matrix $A$ of size $d\times d$, we will denote by $\vert v\vert$ the euclidean norm of $v$ and by $\vert A\vert$ the corresponding operator norm of the matrix $A$. For a random vector $U$ or arbitrary dimension, we write $U\in\L^r$, $r>0$ if $\Vert U\Vert_r=\E^{1/r}\left(\vert U\vert^r\right)<\infty$. 
We will use the following assumptions.

\begin{description}
\item[A1] $\beta$ is a Lipschitz function.
\item[A2] There exist two measurable functions $G:[0,1]\times \left(\R^d\right)^{\N}\rightarrow \R^p$, $H:[0,1]\times \R^{\N}\rightarrow \R$ and a positive constant $C$ such that, if $x_i(u)=G\left(u,\mathcal{F}_i\right)$ and $e_i(u)=H\left(u,\mathcal{F}_i\right)$, 
$$\max_{1\leq i\leq n}\left[\Vert x_i-x_i(i/n)\Vert_2+\Vert x_ie_i-x_i(i/n)e_i(i/n)\Vert_2\right]=O\left(\frac{1}{n}\right),$$
and 
$$\Vert x_i(u)-x_i(v)\Vert_2+\Vert x_i(u)e_i(u)-x_i(v)e_i(v)\Vert_2\leq C\vert u-v\vert.$$ 

\item[A3]
There exist $q>8/7$ and $q_1>16/7$, $q_2>2$ such that $\frac{1}{q}=\frac{1}{q_1}+\frac{1}{q_2}$ and $\Theta_{0,2q_1}(G)+\Theta_{0,q_2}(L)<\infty$, with $L_{n,i}=H_{n,i}G_{n,i}$.
\item[A4]
For each $u\in[0,1]$, we have $\det\E\left[G\left(u,\mathcal{F}_k\right)G\left(u,\mathcal{F}_k\right)'\right]\neq 0$. 
\end{description}

\paragraph{Notes}
\begin{enumerate}
\item
Without additional assumptions, the quantity $\hat{s}_{3,i}$ can be degenerate for finite samples (e.g discrete covariates with the value $0$ in the support of the distribution). 
Assumption {\bf A4} guarantees it is not the case asymptotically. However one can always replace $\hat{s}_{3,i}$ by $\hat{s}_{3,i}+\mu_n I_{p_2}$ where $I_{p_2}$ denotes the identity matrix of size $p_2$ and $(\mu_n)_n$ is a sequence of positive real numbers and decreasing to $0$ at a rate which is sufficiently fast to not modify our asymptotic results (e.g $\mu_n=1/n$). The new matrix is now positive definite and then invertible. Since this modification is 
only technical and complicates our notations, we will not use it.
\item
Assumption {\bf A2} is a local stationarity property. This assumption is satisfied for time-varying autoregressive processes or moving averages with time-varying coefficients. With respect to \citet{ZW}, we directly use the functional dependence coefficient for the triangular arrays $G_{n,i}$, $H_{n,i}$ instead of their local stationary approximations $G\left(i/n,\mathcal{F}_i\right)$ and $H\left(i/n,\mathcal{F}_i\right)$. With this choice, we avoid to state a second result ensuring that the triangular arrays can be replaced in our estimators without changing their asymptotic behavior. Note that when $\left(G_{n,i},H_{n,i}\right)=\left(G\left(i/n,\cdot\right),H\left(i/n,\cdot\right)\right)$, we recover the functional dependence coefficients used in \citet{ZW}. 
\item
The dependence/moment conditions given in assumption {\bf A3} is slightly more restrictive than that of \citet{ZW} which is guaranteed if there exists $q>1$ such that $\Theta_{0,4}(G)+\Theta_{0,2q}(L)<\infty$. For instance, taking $q_1=q_2=2q$, our conditions are satisfied provided $\Theta_{0,4q}(G)+\Theta_{0,2q}(L)<\infty$ for $q>\frac{8}{7}$.

\item
Sufficient conditions for ${\bf A2-A3}$ are given by 
$$\Vert x_i-x_i(u)\Vert_4+\Vert e_i-e_i(u)\Vert_4\leq C\left[\vert u-\frac{i}{n}\vert+\frac{1}{n}\right],$$ 
$$\Vert x_i(u)-x_i(v)\Vert_4+\Vert e_i(u)-e_i(v)\Vert_4\leq C\vert u-v\vert$$
and $\Theta_{0,4q}(G)+\Theta_{0,4q}(H)<\infty$.

\end{enumerate}

For $\ell=1,2$, let us also set $x_{\ell,i}(u)=G_{\ell}\left(u,\mathcal{F}_i\right)$, $e_i(u)=H\left(u,\mathcal{F}_i\right)$
and $y_i(u)=x_{1,i}(u)'\beta_1+x_{2,i}(u)'\beta_2(u)+e_i(u)$. We also set 
$$q_1(u)=\E^{-1}\left[x_{2,i}(u)x_{2,i}(u)'\right]\E\left[x_{2,i}(u)y_i(u)\right],\quad q_2(u)=\E^{-1}\left[x_{2,i}(u)x_{2,i}(u)'\right]\E\left[x_{2,i}(u)x_{1,i}(u)'\right],$$
$$\widetilde{x}_{1,i}(u)=x_{1,i}(u)-q_2(u)'x_{2,i}(u).$$
 
\begin{theo}\label{central}
Assume ${\bf A1-A4}$, $\Theta_{n,2q_1}(G)=O\left(n^{-\nu_1}\right)$ for $\nu_1>\frac{1}{2}-\frac{1}{q_1}$ and $
\Theta_{n,q_2}(L)=O\left(n^{-\nu_2}\right)$ for $\nu_2>\frac{1}{2}-\frac{1}{q_2}$. Then if $b\asymp n^{-c}$ with 
$1/4<c<\min\left\{\frac{3}{4}-q_1^{-1},\frac{3}{4}-q_2^{-1},2-\frac{2}{q\wedge (4/3)},2-\frac{4}{q_1\wedge (8/3)},\frac{1}{2}\right\}$,
we have 
$$\sqrt{n}\left(\hat{\beta}_1-\beta_1\right)\Rightarrow\mathcal{N}\left(0,\Sigma_1^{-1}\Sigma_2\Sigma_i^{-1}\right),$$
with 
$$\Sigma_1=\int_0^1\E\left[\widetilde{x}_{1,i}(u)\widetilde{x}_{1,i}(u)'\right]du,\quad \Sigma_2=\int_0^1\sum_{\ell\in\Z}\c\left[\widetilde{x}_{1,0}(u)e_0(u),\widetilde{x}_{1,\ell}(u)e_{\ell}(u)\right]du.$$
\end{theo}

\paragraph{Notes}
\begin{enumerate}
\item
The optimal bandwidth for estimating Lipschitz time-varying parameters is $b\asymp n^{-1/3}$. 
When $q_1=q_2=2q$ and $q>6/5$, this choice is allowed for applying Theorem \ref{central}.
\item  
It is also possible to use local polynomials of order $1$ for removing the nuisance term $x_{i,2}'\beta_2(i/n)$. 
To this end, let us define for $\ell=0,1,2$,  
$$S_{\ell,i}=\sum_{j=1}^n k_{i,j} x_{2,j}x_{2,j}'\left(\frac{j-i}{nb}\right)^{\ell}.$$
Then we define $\hat{q}_{1,i}$ and $\hat{q}_{2,i}$ as the first components of the two block matrices 
$$\begin{pmatrix} S_{0,i}& S_{1,i}\\S_{1,i}& S_{2,i}\end{pmatrix}^{-1}
\begin{pmatrix} \sum_{j=1}^n k_{i,j} x_{2,j} y_j\\ \sum_{j=1}^n k_{i,j} x_{2,j}y_j\frac{j-i}{nb}\end{pmatrix}\mbox{ and }\begin{pmatrix} S_{0,i}& S_{1,i}\\S_{1,i}& S_{2,i}\end{pmatrix}^{-1}
\begin{pmatrix} \sum_{j=1}^n k_{i,j} x_{2,j} x_{1,j}'\\ \sum_{j=1}^n k_{i,j} x_{2,j}x'_{1,j}\frac{j-i}{nb}\end{pmatrix}$$
respectively.
Moreover, let us denote by $\widetilde{q}_{1,i}$ and $\widetilde{q}_{2,i}$ the first components of the two block matrices
$$\begin{pmatrix} \E S_{0,i}& \E S_{1,i}\\\E S_{1,i}& \E S_{2,i}\end{pmatrix}^{-1}\begin{pmatrix} \sum_{j=1}^n k_{i,j} \E\left[x_{2,j} y_j\right]\\ \sum_{j=1}^n k_{i,j} \E\left[x_{2,j}y_j\right]\frac{j-i}{nb}\end{pmatrix},\quad \begin{pmatrix} \E S_{0,i}& \E S_{1,i}\\\E S_{1,i}& \E S_{2,i}\end{pmatrix}^{-1}\begin{pmatrix} \sum_{j=1}^n k_{i,j} \E\left[x_{2,j} x_{1,j}'\right]\\ \sum_{j=1}^n k_{i,j} \E\left[x_{2,j}x'_{1,j}\right]\frac{j-i}{nb}\end{pmatrix}$$
respectively. Then, under the assumptions of Theorem \ref{central}, one can show that
$$\max_{1\leq i\leq n}\left\vert \widetilde{q}_{1,i}-q_{1,i}\right\vert+\max_{1\leq i\leq n}\left\vert \widetilde{q}_{2,i}-q_{2,i}\right\vert=O\left(b+\frac{1}{n}\right).$$ 
The term $1/n$ corresponds to the deviation with respect to stationarity, i.e the order of the differences $q_{1,i}-q_1(i/n)$ and $q_{2,i}-q_2(i/n)$. 
Moreover we have $\beta_2(i/n)=q_{1,i}-q_{2,i}\beta_1$ and if the second derivative of $\beta_2$ exists and is continuous, we get a smaller bias for $\beta_2$: 
$$\max_{1\leq i\leq n}\left\vert \widetilde{q}_{1,i}-q_{1,i}-\left(\widetilde{q}_{2,i}-q_{2,i}\right)\beta_1\right\vert=O\left(b^2+\frac{1}{n}\right).$$
If we define an estimator of $\beta_1$ with these new expressions for $\hat{q}_{1,i}$ and $\hat{q}_{2,i}$, we obtain the profile least squares estimators studied in \citet{Fan2} for semi-varying coefficients regression models. 
Under the assumptions of Theorem \ref{central}, the estimator $\hat{\beta}_1$ is $\sqrt{n}-$consistent but some assumptions can be relaxed.
Inspection of the proof of Theorem \ref{central} shows that the only change due to this smaller bias concerns the term $N_2$ and the bandwidth condition $nb^4\rightarrow 0$ can be replaced by $nb^6\rightarrow 0$ (the bias is of order $b$ for $\hat{q}_{2,i}$ and of order $b^2$ for $\hat{r}_i$). This leads to the condition $c>1/6$ for the bandwidth. Moreover, on can choose $q>12/11$ and $q_1>24/11$ in assumption {\bf A3} to get a compatible upper bound for $c$ in Theorem \ref{central}.

\end{enumerate}

\subsection{Improving efficiency when the noise process is a martingale difference}
In this subsection, we assume the following martingale difference structure for the noise component $e_i$:
$$\E\left(e_i\big\vert \mathcal{F}_{i-1}\right)=0,\quad \v\left(e_i\big\vert\mathcal{F}_{i-1}\right)=\sigma^2\left(i/n,x_i\right),$$
where $\sigma:[0,1]\times \R^p\rightarrow \R^{*}_+$ is an unknown function. We also assume that $G_{n,i}\left(\mathcal{F}_i\right)=G_{n,i}\left(\mathcal{F}_{i-1}\right)$. 
For simplicity of notations, we set $\sigma_i=\sigma\left(i/n,x_i\right)$. 
Now let $p_i=p\left(i/n,x_i\right)$ a weight function. One can write a second regression model by multiplying the first one by the weights $p_i$, 
\begin{equation}\label{secondmodel}
p_iy_i=p_i x_{1,i}'\beta_1+p_ix_{2,i}'\beta_2(i/n)+p_ie_i,\quad 1\leq i\leq n.
\end{equation}
In the sequel, we assume that the assumptions ${\bf A1-A4}$ are satisfied for this new regression model. 
The asymptotic variance given in Theorem \ref{central} has a simpler form. In particular, 
$$\Sigma_2=\E\int_0^1 p^2_i(u)\sigma^2_i(u)\widetilde{px}_{1,i}(u)\widetilde{px}_{1,i}(u)' du,$$
where $p_i(u)=p\left(u,x_i(u)\right)$, $\sigma_i(u)=\sigma\left(u,x_i(u)\right)$ 
and $px$ denotes the vector $\left(p_ix_i\right)_{1\leq i\leq n}$. 
Now, if we set 
$$q_2(u)=\E^{-1}\left[p_i(u)^2x_{2,i}(u)x_{2,i}(u)'\right]\cdot\E\left[p_i(u)^2 x_{2,i}(u)x_{1,i}(u)'\right],$$
$$q^{*}_2(u)=\E^{-1}\left[\sigma_i(u)^{-2}x_{2,i}(u)x_{2,i}(u)'\right]\cdot\E\left[\sigma_i(u)^{-2} x_{2,i}(u)x_{1,i}(u)'\right]$$
and changing our notations with 
$$\widetilde{x}_{1,i}(u)=x_{1,i}(u)-q_2(u)'x_{2,i}(u),\quad \widetilde{x}_{1,i}^{*}(u)=x_{1,i}(u)-q^{*}_2(u)'x_{2,i}(u),$$
we have 
$$\Sigma_1=\E\int_0^1 p_i(u)^2 \widetilde{x}_{1,i}(u)\widetilde{x}_{1,i}(u)'du=\E\int_0^1p_i(u)^2 \widetilde{x}_{1,i}(u)\widetilde{x}^{*}_{1,i}(u)'du=\E\int_0^1p_i(u)^2\sigma_i(u)\widetilde{x}_{1,i}(u)\frac{\widetilde{x}^{*}_{1,i}(u)'}{\sigma_i(u)}du.$$ 
From the Cauchy-Schwarz inequality, we have for $x,y\in\R^{p_1}$,
$$\left(x'\Sigma_1 y\right)^2\leq x'\Sigma_2 x\cdot y'\Sigma^{*}y,$$
where 
\begin{equation}\label{varopt}
\Sigma^{*}=\E\int_0^1 \frac{\widetilde{x}^{*}_{1,i}(u)\widetilde{x}^{*}_{1,i}(u)'}{\sigma_i(u)^2}du.
\end{equation}
 Now setting $x=\Sigma_1^{-1}z$ and $y=\left(\Sigma^{*}\right)^{-1}z$ for $z\in\R^{p_1}$, we get
$$z'\left(\Sigma^{*}\right)^{-1}z\leq z' \Sigma_1^{-1}\Sigma_2\Sigma_1^{-1}z.$$
Note that $\left(\Sigma^{*}\right)^{-1}$ is the asymptotic variance corresponding to the choice $p_i=\sigma_i^{-1}$. 
This means that weighted least squares can improve efficiency of our estimator, provided to first estimate the conditional variance of the noise. In the subsection $2.4$, we will construct an efficient estimator when the errors are independent but have a time-varying unconditional variance.

\subsection{Comparison with a related method}
\citet{ZW} have proposed another method for the estimation of parameter $\beta_1$. Their method consists in first estimating the full
vector $\beta(\cdot)$ of coefficients using local polynomials and then averaging the first component of this estimator, i.e $\hat{\beta}_1=\int_0^1 \hat{\beta}_1(u)du$. They have shown that under some assumptions, $\hat{\beta}_1$ is root-$n$ consistent. 
Comparing efficiency of both methods seems difficult for a general noise process. 
Nevertheless, one can get interesting properties when the noise is a martingale difference. 
To this end, we keep the notations of the previous subsection. 
Application of the method of \citet{ZW} in the regression model (\ref{secondmodel}) leads to an asymptotic variance $S$ given by (see Theorem $3.1$ of that paper):
$$S=\int_0^1 A S_1(u)^{-1}S_2(u)S_1(u)^{-1}A'du,$$
where 
$$S_1(u)=\E\left[p_i(u)^2x_i(u)x_i(u)'\right],\quad S_2(u)=\E\left[p_i(u)^4\sigma_i(u)^2x_i(u)x_i(u)'\right]$$
and $A$ is a matrix of size $p_1\times p$ and which associates the $p_1$ first components of a vector of $\R^p$. Note that when $p_i=\sigma_i^{-1}$, the asymptotic variance reduces to $\int_0^1A S^{*}(u)^{-1}A' du$ where $S^{*}(u)=\E\left[\frac{x_i(u)x_i(u)'}{\sigma_i(u)^2}\right]$. Using the formula for the inverse of a matrix defined by blocs, one can check that $AS^{*}(u)^{-1}A'=\E^{-1}\left[\frac{h^{*}_i(u)h_i^{*}(u)'}{\sigma_i(u)^2}\right]$. Using Lemma \ref{util1}, 
we have 
\begin{equation}\label{superior}
\int_0^1 AS^{*}(u)^{-1}A'\succeq\left(\Sigma^{*}\right)^{-1}.
\end{equation}
Moreover, the matrix $\int_0^1 AS^{*}(u)^{-1}A'$ is optimal for the estimator of \citet{ZW}. Indeed we have from the Cauchy-Schwarz 
inequality 
$$\left(x'S_1(u)x\right)^2\leq x'S_2(u)x\cdot y'S^{*}(u)y$$
and the result follows by tacking $x=S_1(u)^{-1}z$ and $y=S^{*}(u)^{-1}z$.  
Then we obtain the following important conclusion for a martingale difference noise: removing the nonparametric component with partial regressions can lead to a more efficient estimator than averaging directly the nonparametric estimator.
In the next section, we show that an AR process with a conditionally Gaussian noise and some strongly exogenous covariates, the matrix $\left(\Sigma^{*}\right)^{-1}$ is a lower bound in semiparametric estimation.

\subsection{Weighted least squares and asymptotic efficiency improvement}\label{plugin}
In this subsection, we consider the simple case $e_i=\sigma(i/n)f\left(\xi_i\right)$ with $\E f\left(\xi_i\right)=0$ and $\v f\left(\xi_i\right)=1$. In this case, we first derive an estimator of $\sigma^2(i/n)$. A natural candidate is 
$$\hat{\sigma}^2(i/n)=\sum_{j=1}^n k_{i,j}\left[y_j-x_j'\hat{\beta}(i/n)\right]^2,$$
where $\hat{\beta}$ is an estimate of $\beta$. We will use the simple Nadaraya-Watson estimator of $\beta$, i.e 
$$\hat{\beta}(i/n)=\left(\sum_{j=1}^n k_{i,j}x_jx_j'\right)^{-1}\sum_{j=1}^nk_{i,j}x_jy_j.$$
Another solution is to take the semiparametric estimator $\hat{\beta}_1$ of $\beta_1$ which is $\sqrt{n}-$consistent under the assumptions of Theorem \ref{central} and to define an estimator of $\beta_2(i/n)$ as follows:
$$\hat{\beta}_2(i/n)=\arg\min_{\alpha\in\R^{p_2}}\sum_{j=1}^n k_{i,j}\left[y_i-x_{i,1}'\hat{\beta}_1-x_{i,2}'\alpha\right]^2.$$
But this second choice will not improve our asymptotic results and for simplicity we only consider the Nadaraya-Watson estimator for the plug-in.
Our plug-in estimator is defined as follows. For simplicity we will simply note $\sigma_i$ (resp. $\hat{\sigma}_i$) for $\sigma(i/n)$ (resp. $\hat{\sigma}(i/n)$). 
We set 
$$\hat{s}^{*}_{1,i}=\sum_{j=1}^nk_{i,j}\frac{x_{2,j}y_j}{\hat{\sigma}_j^2},\quad \hat{s}^{*}_{2,i}=\sum_{j=1}^nk_{i,j}\frac{x_{2,j}x_{1,j}'}{\hat{\sigma}_j^2},\quad \hat{s}^{*}_{3,i}=\sum_{j=1}^nk_{i,j}\frac{x_{2,j}x_{2,j}'}{\hat{\sigma}_j^2},$$
$$\hat{q}^{*}_{1,i}=\left(\hat{s}^{*}_{3,i}\right)^{-1}\hat{s}^{*}_{1,i},\quad \hat{q}^{*}_{2,i}=\left(\hat{s}^{*}_{3,i}\right)^{-1}\hat{s}^{*}_{2,i}$$
and 
$$\hat{y}^{*}_i=y_i-x_{2,i}'\hat{q}^{*}_{1,i},\quad \hat{x}^{*}_{1,i}=x_{1,i}-\left(\hat{q}^{*}_{2,i}\right)'x_{2,i}.$$
The optimal estimator of $\beta_1$ is defined by 
$$\hat{\beta}_1^{*}=\left(\sum_{i=1}^n\frac{\hat{x}^{*}_{1,i}\left(\hat{x}^{*}_{1,i}\right)'}{\hat{\sigma}_i^2}\right)^{-1}\sum_{i=1}^n \frac{\hat{x}^{*}_{1,i}\hat{y}^{*}_i}{\hat{\sigma}_i^2}.$$

We have the following result.

\begin{theo}\label{central2}
Assume that all the assumptions {\bf A1-A4} hold true with $q_2>4$ and $\sigma$ is positive and Lipschitz over $[0,1]$. 
Then
$$\sqrt{n}\left(\hat{\beta}^{*}_1-\beta_1\right)\Rightarrow \mathcal{N}\left(0,\left(\Sigma^{*}\right)^{-1}\right),$$
where $\Sigma^{*}$ is defined in (\ref{varopt}).
\end{theo}

\section{Asymptotic semiparametric efficiency for an autoregressive process with exogenous covariates}

\subsection{An autoregressive process with covariates}

We consider a process $\left(z_i\right)_{1\leq i\leq n}$ such that $z_i=G^{(2)}_{n,i}\left(\mathcal{F}_i\right)$ and an error process 
$(e_i)_{1\leq i\leq n}$ defined as in Section $2$.
 We define a process $(y_i)_{1\leq i\leq n}$ recursively by 
$$y_i=\sum_{j=1}^q a_j(i/n)y_{i-j}+z_i'\gamma(i/n)+e_i$$
and we set $x_i=\left(y_{i-1},\ldots,y_{i-q},z_i'\right)'$ and $p=q+d$. Finally, we set $\beta(u)=\left(a_1(u),\ldots,a_q(u),\gamma(u)'\right)'$ for the parameter of interest. 
The following assumptions will be needed.
\begin{description}
\item[B1]
The function $u\mapsto \beta(u)$ is Lipschitz. Moreover, for each $u\in[0,1]$, the roots of the polynomial $\mathcal{P}_u(z)=1-\sum_{j=1}^q a_j(u)z^j$ are outside the unit disc.
\item[B2]
There exists a measurable function $M:[0,1]\times E^{\N}\rightarrow \R$ such that if $z_i(u)=M\left(u,\mathcal{F}_i\right)$,
$$\max_{1\leq i\leq n}\left[\Vert z_i-z_i(i/n)\Vert_4+\Vert e_i-e_i(i/n)\Vert_4\right]=O\left(\frac{1}{n}\right),$$
and
$$\Vert z_i(u)-z_i(v)\Vert_4+\Vert e_i(u)-e_i(v)\Vert_4\leq C\vert u-v\vert.$$
\item[B3]
There exist $q>\frac{8}{7}$ such that $\Theta_{0,4q}\left(M\right)+\Theta_{0,4q}(H)<\infty$.
\end{description}
Under the previous assumptions, we show that $\left(y_i\right)_{1\leq i\leq n}$ can be approximated locally by the stationary process
$$y_i(u)=\sum_{j=1}^qa_j(u)y_{i-j}(u)+z_i(u)'\gamma(u)+e_i(u),\quad i\in\Z,u\in [0,1].$$

For initialization, we assume that for $i\leq 0$, $y_i=y_i(0)$. Note that for integers $n\geq 1$ and $1\leq i\leq n$, 
there exists an application $N_{n,i}:\left(\R^d\right)^{\N}\rightarrow \R$ such that $y_i=N_{n,i}\left(\mathcal{F}_i\right)$.
Moreover, there exists also an application $N:[0,1]\times \left(\R^d\right)^{\N}\rightarrow \R$ such that $y_i(u)=N\left(u,\mathcal{F}_i\right)$.

\begin{prop}\label{approximationAR}
Under the assumptions {\bf B1-B3}, there exists a constant $C>0$ such that for all $(u,v,i)\in [0,1]^2\times \{1,2,\ldots,n\}$,
$$\Vert y_i-y_i(u)\Vert_4\leq C\left[\vert u-\frac{i}{n}\vert+\frac{1}{n}\right],\quad \Vert y_i(u)-y_i(v)\Vert_4\leq C\vert u-v\vert.$$
Moreover, we have $\Theta_{0,4q}(N)<\infty$.
\end{prop}

\paragraph{Note.} Assumptions {\bf B1-B3} entail assumptions {\bf A1-A3}. There is no general arguments for checking {\bf A4}. However, if {\bf A4} is not valid then we have $\sum_{j=1}^q\lambda_jy_{i-j}(u)+\mu'z_i(u)=0$ a.s for some $u\in [0,1]$ and some constants $\lambda_j$ and $\mu$. This type of equality cannot hold in general. For instance, assumption {\bf A4} is automatic for a non degenerate noise such that the $e_i'$s are independent and independent from the $z_i'$s. 

\paragraph{Proof of Proposition \ref{approximationAR}}
Let us introduce the following notations. For $i\leq n$, we set 
$$Y_i=\left(y_i,y_{i-1},\ldots,y_{i-q+1}\right)',\quad  B_i=\left(z_i'\gamma(i/n)+e_i,{\bf 0}_{1,q-1}'\right)',$$
with the convention that $y_i=y_i(0)$ and $\left(\gamma(i/n)',e_i,z_i'\right)=\left(\gamma(0)',e_i(0),z_i(0)\right)$ when $i\leq 0$. 
Moreover we denote by $A(i/n)$ the companion matrix associated to the polynomial $\mathcal{P}_{i/n}$. 
Then, for $q+1\leq i\leq n$, we have $Y_i=A(i/n)Y_{i-1}+B_i$. From the assumption {\bf B2}, for each $u\in [0,1]$, the spectral radius 
of the matrix $A(u)$ is smaller than $1$. Then, one can consider
$$m_u=\inf\left\{k\geq 1: \left\vert A(u)^k\right\vert <1\right\}<\infty.$$
Now, let $\mathcal{O}_u=\left\{v\in [0,1]: \left\vert A(v)^{m_u}\right\vert<1\right\}$. The set $\mathcal{O}_u$ is an open subset of $[0,1]$. From the compactness of the unit interval, we have $[0,1]=\cup_{i=1}^h\mathcal{O}_{u_i}$ for suitable points $u_1,\ldots,u_h$ in $[0,1]$. Then if $m$ denotes the lowest common multiple of the integers $m_{u_1},\ldots,m_{u_h}$, we have $\max_{u\in [0,1]}\left\vert A(u)^m\right\vert<1$. From the uniform continuity of the application $(u_1,\ldots,u_h)\mapsto A(u_1)\cdots A(u_h)$, we have, if $n$ is large enough, $$\rho=\max_{m-1\leq k\leq n}\left\vert A\left(\frac{k}{n}\right)A\left(\frac{k-1}{n}\right)\cdots A\left(\frac{k-m+1}{n}\right)\right\vert<1.$$
Moreover, we have $\sup_{u\in [0,1]}\Vert Y_i(u)\Vert_4<\infty$.  
From our assumptions, we deduce that there exists a positive real number $D$ such that, 
\begin{eqnarray*}
\Vert Y_i-Y_i(u)\Vert_4&\leq& \rho\Vert Y_{i-m}-Y_{i-m}(u)\Vert_4+D\left[\left\vert\frac{i}{n}-u\right\vert+\frac{1}{n}\right]\\
&\leq & D\sum_{h=0}^{\infty}\rho^h\left[\left\vert u-\frac{i-h\ell}{n}\right\vert+\frac{1}{n}\right],
\end{eqnarray*}
which leads to the first part of the proposition, the bound for $\Vert Y_i(u)-Y_i(v)\Vert_4$ can be obtained in the same way. 
For the second part, we use the expansion
$$Y_i=B_i+\sum_{h=1}^{\infty}A\left(\frac{i}{n}\right)\cdots A\left(\frac{i-h+1}{n}\right)B_{i-h}.$$
Then we obtain $\Theta_{0,4q}(N)<\infty$ by using the contraction coefficient $\rho$ and the assumption $\Theta_{0,4q}(M)+\Theta_{0,4q}(H)<\infty$.$\square$

\subsection{Gaussian inputs and efficient information matrix}
In this subsection, we assume that $e_i=\sigma(i/n,z_i)f\left(\xi_i\right)$ with $f(\xi_i)\sim\mathcal{N}(0,1)$ and the two processes 
$(z_i)_{1\leq i\leq n}$ and $\left(f(\xi_i)\right)_{i\in\Z}$ are independent. This independence property occurs for example when 
$\xi_i=\left(\xi_i^{(1)},\xi_i^{(2)}\right)$, with $\xi^{(1)}$ and $\xi^{(2)}$ independent, $f(\xi_i)=\xi_i^{(1)}$ and $z_i$ only depending 
on $\left(\xi_j^{(2)}\right)_{j\leq i}$. This type of condition on the covariates is called strong exogeneity (sometimes Granger or Sims causality in the literature) and is very useful to get a simple expression of the conditional likelihood function.
We will assume that the conditional standard deviation $\sigma_i=\sigma\left(i/n,z_i\right)$ satisfies $\inf_{n\geq 1,1\leq i\leq n}\sigma_i\geq \sigma_{-}>0$ a.s. for a constant $\sigma_{-}$. One could also allow $\sigma_i$ to depend on past values of $y_i$ but 
it is necessary to introduce additional regularity conditions to define the model. This is why we prefer to consider conditional heteroscedasticity only with respect to the covariates.

We assume that the parameter $\beta(u)$ writes as $\left(\beta_1,\beta_2(u)\right)$ with $\beta_1\in \R^{p_1}$ is constant. 
The coordinates of $\beta_1$ can contain some of the $a_j'$s or some of the coordinates of $\gamma$. 
For instance, $\beta_1$ can contains all the parameters except a time-varying trend $\beta_2$. One can also consider the case of constant lag coefficients $a_j$ and a covariate $z_i$ with time-varying coefficients.

We will first prove a LAN expansion of the likelihood ratio conditionally to some initial values $y_0,y_1,\ldots,y_{-q+1}$ not depending on $n$. Let $\mu_n$ be the probability distribution of the vector $\left(z_1',z_2',\ldots,z_n'\right)'$. Note that $\mu_n$ does not depend on a parameter of interest in this problem. We also denote by $L_{n,\beta}\left(\cdot\vert y_0^{-},z_{1:n}\right)$ the density of $(y_1,y_2,\ldots,y_n)$ conditionally to $z_i$ and $y_j$ for $1\leq i\leq n$ and $-q+1\leq j\leq 0$. Using the exogeneity property of the covariates, we have the expression 
$$L_{n,\beta}\left(y_{1:n}\vert y_0^{-},z_{1:n}\right)=\prod_{i=1}^n \left(2\pi\sigma_i^2\right)^{-1/2}\exp\left(-\frac{\left\vert y_i-x_i'\beta(i/n)\right\vert^2}{2\sigma_i^2}\right).$$  
Then we set 
$$d P_{n,\beta}\left(y_{1:n},z_{1:n}\vert y_0^{-}\right)=L_{n,\beta}\left(y_{1:n}\vert y_0^{-},z_{1:n}\right)d\mu_n\left(z_{1:n}\right)d\lambda_n\left(y_{1:n}\right),$$
where $\lambda_n$ denotes the Lebesgue measure on $\R^n$. 
We introduce the following local parametrization. Let us denote by $\mathbb{H}$ the Hilbert space $\R^{p_1}\times\mathcal\L^{p_1}([0,1])$ endowed with the scalar product 
$$\left(h\vert g\right)_{\mathbb{H}}=\int_0^1 h(u)'M(u)h(u)du,\quad M(u)=\E\left[\frac{x_0(u)x_0(u)'}{\sigma^2\left(u,x_0(u)\right)}\right].$$
We set $h(u)$, $g(u)$ for simplicity of notations. It is intended that the $p_1-$first components of $h(u)$ and $g(u)$ do not depend on $u$.
Note that, in a parametric problem, $M(u)$ corresponds to the Fisher information matrix for the estimation of the parameter $\beta(u)$ in the model
$$y_i(u)=x_i'(u)\beta(u)+\sigma\left(u,z_i(u)\right)f\left(\xi_i\right).$$

We also set $\mathcal{H}=\R^{p_1}\times \mathcal{L}^{p_2}$ where $\mathcal{L}$ denotes the set of Lipschitz functions defined over $[0,1]$.
Then if $\beta,h\in \mathcal{L}\times\mathcal{L}$, the log likelihood ratio obtained after shifting parameter $\beta$ is given by 
$$\log\frac{dP_{n,\beta+h/\sqrt{n}}}{dP_{n,\beta}}\left(y_{1:n},z_{1:n}\vert y_0^{-}\right)=\sum_{i=1}^n \left[\frac{\left\vert y_i-x_i'\beta\right\vert^2}{2\sigma_i^2}-\frac{\left\vert y_i-x_i'(\beta+h)\right\vert^2}{2\sigma_i^2}\right].$$   
Let us first derive a LAN expansion for this likelihood ratio.

\begin{prop}\label{LAN}
Under the assumptions ${\bf B1-B3}$, we have
$$\log\frac{dP_{n,\beta+h/\sqrt{n}}}{dP_{n,\beta}}\left(y_{1:n},z_{1:n}\vert y_0^{-}\right)=\Delta_{n,\beta,h}-\frac{1}{2}\Vert h\Vert^2_{\mathbb{H}}+o_{\P}(1),$$
where 
$$\Delta_{n,\beta,h}:=\frac{1}{\sqrt{n}}\sum_{i=1}^n\frac{f\left(\xi_i\right)x_i'h(i/n)}{\sigma_i}\Rightarrow \mathcal{N}\left(0,\Vert h\Vert^2_{\mathbb{H}}\right).$$
\end{prop}

\paragraph{Proof of Proposition \ref{LAN}}
We have 
$$\log\frac{dP_{n,\beta+h/\sqrt{n}}}{dP_{n,\beta}}\left(y_{1:n},z_{1:n}\vert y_0^{-}\right)=\Delta_{n,\beta,h}-\frac{1}{2}\Delta'_{n,h},$$
where 
$$\Delta'_{n,h}=\frac{1}{n}\sum_{i=1}^n \frac{h'(i/n)x_ix_i'h(i/n)}{\sigma_i^2}.$$
First, one can note that the initial values $y_0,y_{-1},\ldots,y_{-q+1}$ are forgotten exponential fast. More precisely if $y_i$ is initialized with the random variables $y_0(0),y_{-1}(0),\ldots,y_{-q+1}(0)$, the $\L^4-$norm of the difference is bounded by $C\rho^i$ with $\rho\in (0,1)$. Hence, one can assume that $(y_i)_{1\leq i\leq n}$ is initialized with $\left(y_i(0)\right)_{-q+1\leq i\leq 0}$. 
Using the assumptions ${\bf B1-B3}$ and Proposition \ref{approximationAR}, there exists a positive constant $C$ such that for all $n\geq 1$, $1\leq i\leq n$ and $u\in [0,1]$,
$$\left\vert \E\left(\sigma_i^{-2}x_ix_i'\right)-\E\left(\sigma_i(u)^{-2}x_i(u)x_i(u)'\right)\right\vert\leq C\left[\vert u-i/n\vert+1/n\right].$$
Moreover, the application $u\mapsto \E\left(\sigma_i(u)^{-2}x_i(u)x_i(u)'\right)$ is Lipschitz and  
$\max_{1\leq i\leq n} \E \sigma_i^{-4}\vert x_i\vert^4=O(1)$. 
Then the convergence of $\Delta_{n,\beta,h}$ follows from the central limit theorem for triangular arrays of martingale differences. 
Moreover, using Lemma $A1$ of \citet{ZW} with $q=2$, we have
$$\Delta'_{n,h}-\frac{1}{n}\sum_{i=1}^nh(i/n)'\E\left[\sigma_i^{-2}x_ix_i'\right]h(i/n)=O_{\P}\left(1/\sqrt{n}\right).$$
Then we deduce from the properties stated above that $\Delta'_{n,h}\rightarrow \Vert h\Vert^2_{\mathbb{H}}$. This completes the proof of Proposition \ref{LAN}.$\square$

From this LAN expansion, we deduce the following result.

\begin{cor}\label{opti}
The limit distribution of any regular estimator of parameter $\beta_1$ equals the distribution of $L_1+L_2$ of independent random vectors such that 
$$L_1\sim\mathcal{N}_{p_1}\left(0,\left(\Sigma^{*}\right)^{-1}\right),\quad \Sigma^{*}=\int_0^1\E\left(\sigma_i(u)^{-2}\widetilde{x}^{*}_{1,i}(u)\widetilde{x}^{*}_{1,i}(u)'\right)du,$$
with 
$$\widetilde{x}^{*}_{1,i}(u)=x_{1,i}(u)-q_2^{*}(u)'x_{2,i}(u),\quad q_2^{*}(u)=\E^{-1}\left[\sigma_i(u)^{-2}x_{2,i}(u)x_{2,i}(u)'\right]\E\left[ \sigma_i(u)^{-2}x_{2,i}(u)x_{1,i}(u)'\right].$$
\end{cor}

\paragraph{Note.} The matrix $\Sigma^{*}$ plays the rule of an efficient information matrix. Its definition involves the whole family of stationary approximations $\left\{\left(y_i(u),z_i(u),e_i(u)\right):u\in [0,1]\right\}$ of our triangular array.

\paragraph{Proof of Corollary \ref{opti}}
We use the general convolution theorem given in \citet{VWW}, theorem $3.11.2$, which is compatible with triangular arrays. 
To this end, we consider the sequence of parameters $\left\{\kappa_n(h): h\in H\right\}$ defined by $\kappa_n(h)=\beta_1+h_1/\sqrt{n}$ where $h_1$ denotes the vector of constant parameters. This sequence of parameters is regular, we have 
$$\sqrt{n}\left(\kappa_n(h)-\kappa_n(0)\right)=h_1:=\dot{\kappa}(h).$$
The adjoint operator $\dot{\kappa}^{*}:\R^{p_1}\rightarrow \H$ of the projection operator $\dot{\kappa}$ is defined by the equalities, 
$$\dot{\kappa}(h)'v=\left(h,\dot{\kappa}^{*}v\right)_{\H},\quad v\in\R^{p_1},h\in H.$$
First, we replace the two random vectors $x_{1,i}(u)/\sigma_i(u)$ and $x_{2,i}(u)/\sigma_i(u)$ by orthogonal vectors in $\L^2$. We have 
$$\frac{x_{1,i}(u)}{\sigma_i(u)}=\frac{\widetilde{x}_{1,i}^{*}(u)}{\sigma_i(u)}+q^{*}_2(u)'\frac{x_{2,i}(u)}{\sigma_i(u)}.$$
Then we have 
$$M(u)=\begin{pmatrix} I_{p_1} & q^{*}_2(u)'\\{\bf \Large 0}_{p_2,p_1}& I_{p_2}\end{pmatrix}\cdot\begin{pmatrix} \E\left[\frac{\widetilde{x}_{1,i}^{*}(u)\widetilde{x}_{1,i}^{*}(u)'}{\sigma_i(u)^2}\right]& {\bf \Large 0}_{p_1,p_2}\\ {\bf \Large 0}_{p_2,p_1}& \E\left[\frac{x_{2,i}(u)x_{2,i}(u)'}{\sigma_i(u)^2}\right]\end{pmatrix}\cdot\begin{pmatrix} I_{p_1}& {\bf \Large 0}_{p_1,p_2}\\q^{*}_2(u)& I_{p_2}\end{pmatrix}.$$
Then if $\dot{\kappa}^{*}v=(g_1',g_2(u)')'$, we get 
$$h_1'v=h_1'\Sigma g_1+\int_0^1\left(h_1'q^{*}_2(u)'+h_2(u)'\right)\E\left[\frac{x_{2,i}(u)x_{2,i}(u)'}{\sigma_i(u)^2}\right]\left(q^{*}_2(u)g_1+g_2(u)\right)du.$$
Then we get $g_1=\Sigma^{-1}v$ and $g_2(u)=-q^{*}_2(u)\Sigma^{-1}v$. 
We deduce that 
$$\Vert \dot{\kappa}^{*}v\Vert_{\H}=v'\Sigma^{-1}v$$
and the convolution theorem gives the result, from the LAN expansion given in Proposition \ref{LAN}.$\square$

\paragraph{Note.} Under the assumptions of Theorem $2$ and when $\sigma_i=\sigma(i/n)$, we deduce that the estimator $\hat{\beta}^{*}_1$ is asymptotically efficient
in a semiparametric sense. Rigorously, it is necessary to show that this estimator is regular in the sense of \citet{VWW}. This step does not present any difficulty, provided to follow carefully the proof of Theorem $2$. Details are omitted.

\section{Proofs of Theorem \ref{central} and Theorem \ref{central2}} 

\subsection{Additional results for the proof of Theorem \ref{central}}
In the sequel, we set $\widetilde{q}_{\ell,i}=\E^{-1}\left(\hat{s}_{3,i}\right)\E\left(\hat{s}_{\ell,i}\right)$ for $\ell=1,2$.
\begin{lem}\label{approxquotient}
Assume that the assumptions of Theorem \ref{central} hold true. Setting $m_i=\E\hat{s}_{3,i}$, we have
\begin{enumerate}
\item
$\max_{1\leq i\leq n}\Vert m_i^{-1}\Vert=O(1)$,  $\max_{1\leq i\leq n}\Vert \hat{s}_{3,i}^{-1}\Vert=O_{\P}(1)$,  
\item
for $\ell=1,2$,  
$$\max_{1\leq i\leq n}\left\vert \hat{q}_{\ell,i}-\widetilde{q}_{\ell,i}-m_i^{-1}\left[\hat{s}_{\ell,1}-\E\hat{s}_{\ell,j}\right]+
m_i^{-1}\left[\hat{s}_{3,i}-\E\hat{s}_{3,i}\right]m_i^{-1}\E\hat{s}_{\ell,i}\right\vert=o_{\P}\left(\frac{1}{\sqrt{n}}\right).$$
\item
for $\ell=1,2$, $\max_{1\leq i\leq n}\left\vert \hat{q}_{\ell,i}-\widetilde{q}_{\ell,i}\right\vert=o_{\P}\left(n^{-1/4}\right).$
\item
For $\ell=1,2$, we have 
$$\max_{1\leq i\leq n}\left\vert \widetilde{q}_{\ell,i}-q_{\ell,i}\right\vert=o\left(n^{-1/4}\right),\quad 
\max_{1\leq i\leq n}\left\vert q_{\ell,i}-q_{\ell}(i/n)\right\vert=O(1/n).$$
\end{enumerate}
\end{lem}

\paragraph{Proof of Lemma \ref{approxquotient}}
\begin{enumerate}
\item
Using Lemma A1 in \citet{ZW}, there exists $C>0$ such that  
$$\Vert\max_{1\leq i\leq n}\left\vert\hat{s}_{3,i}-\E\hat{s}_{3,i}\right\vert\Vert_2\leq C \sqrt{n} a_n \Theta_{0,4}(G),$$
where 
$$a_n=\max_{1\leq i\leq n}\left[k_{i,1}+\sum_{j=1}^{n-1}\vert k_{i,j+1}-k_{i,j}\vert\right]
= O\left(\frac{1}{nb}\right).$$

Using our assumptions, we deduce that 
\begin{equation}\label{etap1}
\max_{1\leq i\leq n}\left\vert \hat{s}_{3,i}-\E\hat{s}_{3,i}\right\vert=o_{\P}(1).
\end{equation}
Moreover, setting $G=(G_1,G_2)$ and $g_2(u)=\E G_2\left(u,\mathcal{F}_i\right)G_2\left(u,\mathcal{F}_i\right)'$, we have from assumption {\bf A2},
\begin{equation}\label{etap2}
\max_{1\leq i\leq n}\left\vert \E\hat{s}_{3,i}-g_2(i/n)\right\vert=O(b+1/n).
\end{equation}
From (\ref{etap1}) and (\ref{etap2}), we conclude that $\max_{1\leq i\leq n}\left\vert \hat{s}_{3,i}-g_2(i/n)\right\vert=o_{\P}(1)$. 
Then the result follows from assumptions {\bf A2} and {\bf A4} which guarantee that $\max_{1\leq i\leq n}\left\vert g_2(i/n)^{-1}\right\vert=O(1)$.
\item
One can check that the quantity inside the norm is equal to 
\begin{eqnarray*}
d_i&=&\hat{s}_{3,i}^{-1}\left(\hat{s}_{3,i}-\E\hat{s}_{3,i}\right)m_i^{-1}\left(\hat{s}_{3,i}-\E\hat{s}_{3,i}\right)m_i^{-1}\hat{s}_{\ell,i}\\&+& m_i^{-1}\left(\E\hat{s}_{3,i}-\hat{s}_{3,i}\right)m_i^{-1}\left(\hat{s}_{\ell,i}-\E\hat{s}_{\ell,i}\right).
\end{eqnarray*} 
Using Lemma \ref{uniforme} and the previous point of this lemma, we get 
$$\max_{1\leq i\leq n}\vert d_i\vert=O_{\P}\left(\frac{v_n^2}{(nb)^2}\right),\mbox{ where } v_n=n^{1/q_1}+n^{1/q_2}+\sqrt{nb \log n}$$
and the result follows from our bandwidth conditions.

\item
The result follows from the previous point, Lemma \ref{uniforme} and our bandwidth conditions.

\item
As in the proof of point $1.$, for $\ell=1,2,3$, we have from assumption {\bf A2}  
$$\max_{1\leq i\leq n}\left\vert \E\hat{s}_{\ell,i}-\sum_{j=1}^n k_{i,j}s_{\ell}(j/n)\right\vert=O\left(1/n\right),\quad \max_{1\leq i\leq n}\left\vert \sum_{j=1}^nk_{i,j}s_{\ell}(j/n)-s_{\ell}(i/n)\right\vert=O(b).$$
Then the result follows from assumption {\bf A4} and the condition $b=o\left(n^{-1/4}\right)$.$\square$
  
\end{enumerate}

\subsection{Proof of Theorem \ref{central}}
We set $D_n=\sum_{i=1}^n \hat{x}_{1,i}\hat{x}_{1,i}'$.
Using the decomposition 
$$\hat{y}_i=\hat{x}_{1,i}'\beta_1+x_{2,i}'\hat{r}_i+e_i,\quad \hat{r}_i=\left(\hat{q}_{2,i}-q_{2,i}\right)\beta_1-\left(\hat{q}_{1,i}-q_{1,i}\right),$$
we get 
$$\hat{\beta}_1-\beta_1=D_n^{-1}\sum_{i=1}^n \hat{x}_{1,i}\left[x_{2,i}'\hat{r}_i+e_i\right].$$
Using now the decomposition $\hat{x}_{1,i}=\widetilde{x}_{1,i}-\left(\hat{q}_{2,i}-q_{2,i}\right)'x_{2,i}$, we have
$$\hat{\beta}_1-\beta_1=D_n^{-1}\left[N_1-N_2+N_3-N_4\right],$$
where 
$$N_1=\sum_{i=1}^n \widetilde{x}_{1,i} x_{2,i}'\hat{r}_i,\quad N_2=\sum_{i=1}^n \left(\hat{q}_{2,i}-q_{2,i}\right)'x_{2,i}x_{2,i}'\hat{r}_i,$$
$$N_3=\sum_{i=1}^n \widetilde{x}_{1,i}e_i,\quad N_4=\sum_{i=1}^n\left(\hat{q}_{2,i}-q_{2,i}\right)'x_{2,i}e_i.$$
The proof of Theorem \ref{central} will follow from the following points.
\begin{itemize}
\item
We first show that the leading term $\frac{N_3}{\sqrt{n}}=\frac{1}{\sqrt{n}}\sum_{i=1}^n\widetilde{x}_{1,i}e_i$ converges in distribution. 
Using assumption {\bf A2} and Lemma \ref{approxquotient} (4.), we have 
$$\frac{1}{\sqrt{n}}\sum_{i=1}^n \Vert \widetilde{x}_{1,i}e_i-\widetilde{x}_{1,i}(i/n)e_i(i/n)\Vert_2=O\left(\frac{1}{\sqrt{n}}\right).$$
As in the proof of Theorem $1$ given in \citet{ZW}, using $m-$dependent approximations, it is possible to show that 
$$\frac{1}{\sqrt{n}}\sum_{i=1}^n\widetilde{x}_{1,i}(i/n)e_i(i/n)\Rightarrow \mathcal{N}\left(0,\Sigma_2\right).$$
\item
Using Lemma \ref{approxquotient}, we have 
$$\frac{1}{\sqrt{n}}\left\vert N_2\right\vert\leq C\frac{1}{\sqrt{n}}\sum_{i=1}^n \left\vert x_{2,i}x_{2,i}'\right\vert\cdot\max_{\ell=1,2}\max_{1\leq i\leq n}\left\vert \hat{q}_{\ell,i}-q_{\ell,i}\right\vert^2=o_{\P}(1).$$ 

\item
Next, we show that $\frac{N_4}{\sqrt{n}}=o_{\P}(1)$. With the same arguments, one can also show that $\frac{N_1}{\sqrt{n}}=o_{\P}(1)$. 
Using Lemma \ref{approxquotient}, it is only necessary to prove that $\frac{1}{\sqrt{n}}\sum_{i=1}^n A_i'x_{2,i}e_i=o_{\P}(1)$ when 
$$A_i=\widetilde{q}_{2,i}-q_{2,i},\quad A_i=m_i^{-1}\left(\hat{s}_{2,i}-\E\hat{s}_{2,i}\right),\quad A_i=m_i^{-1}\left(\hat{s}_{3,i}-\E\hat{s}_{3,i}\right)m_i^{-1}\E\hat{s}_{2,i}.$$
The first case follows from Lemma \ref{approxquotient}, Lemma $A.1$ in \citet{ZW} and our bandwidth conditions.  
The second and third case follows from Lemma \ref{smooth1}, applied with $r_1=q_1$, $r_2=q_2$ and using our bandwidth conditions.
One can also note that for the proof of $N_1=o\left(\sqrt{n}\right)$, we have to show that $\frac{1}{\sqrt{n}}\sum_{i,j=1}^n k_{i,j}z_i w_j=o_{\P}(1)$ when $z_i$ and $w_j$ are elements of the matrices $x_ix_i'$ and $x_jx_j'$ respectively. One can apply Lemma \ref{smooth1}
with $r_1=r_2=q_1$ and the result follows from our bandwidth conditions.

\item
To end the proof of Theorem \ref{central}, it remains to show that $\frac{D_n}{n}=\Sigma_1+o_{\P}(1)$.
We first note that 
\begin{eqnarray*}
D_n&=&\sum_{i=1}^n \widetilde{x}_{1,i}\widetilde{x}_{1,i}'+\sum_{i=1}^n\left(\hat{q}_{2,i}-q_{2,i}\right)x_{2,i}x_{2,i}'\left(\hat{q}_{2,i}-q_{2,i}\right)'\\
&+& \sum_{i=1}^n\left(\hat{q}_{2,i}-q_{2,i}\right)' x_{2,i}\widetilde{x}_{1,i}'+\sum_{i=1}^n \widetilde{x}_{1,i}x_{2,i}'\left(\hat{q}_{2,i}-q_{2,i}\right).
\end{eqnarray*}
It has already been shown in the previous points that 
$$\frac{D_n}{n}=\frac{1}{n}\sum_{i=1}^n \widetilde{x}_{1,i}\widetilde{x}_{1,i}'+o_{\P}(1).$$
Moreover, using the moment bound given by Lemma $A.1$ in \citet{ZW}, we have 
$$\frac{1}{n}\sum_{i=1}^n \left[\widetilde{x}_{1,i}\widetilde{x}_{1,i}'-\E\left(\widetilde{x}_{1,i}\widetilde{x}_{1,i}'\right)\right]=o_{\P}(1).$$
Then using assumption {\bf A2} and the convergence of Riemann sums, we get
\begin{eqnarray*}
\frac{1}{n}D_n&=&\frac{1}{n}\sum_{i=1}^n \E\left[\widetilde{x}_{1,0}(i/n)\widetilde{x}_{1,0}(i/n)'\right]+o_{\P}(1)\\
&=& \Sigma_1+o_{\P}(1).\square
\end{eqnarray*}
\end{itemize}

\subsection{Additional results for the proof of Theorem \ref{central2}}
\begin{lem}\label{plugin}
Assume that the assumptions of Theorem \ref{central2} hold true. Then the following properties are valid.
\begin{enumerate}
\item
$\max_{1\leq i\leq n}\left\vert \hat{\beta}(i/n)-\beta(i/n)\right\vert=O_{\P}\left(b+\frac{\sqrt{\log(n)}}{\sqrt{nb}}+\frac{1}{n^{1-\frac{1}{q_2}}b}\right)=o_{\P}\left(n^{-1/4}\right).$
\item
$\max_{1\leq i\leq n}\sum_{j=1}^n k_{i,j}\vert x_j\vert^2=O_{\P}(1)$ and $\max_{1\leq i\leq n}\sum_{j=1}^n k_{i,j}\vert x_je_j\vert=O_{\P}(1)$.
\item
$\max_{1\leq i\leq n}\left\vert \hat{\sigma}_i^2-\sum_{j=1}^n k_{i,j}e_j^2\right\vert=o_{\P}\left(\frac{1}{\sqrt{n}}\right)$ and 
$\max_{1\leq i\leq n}\left\vert \sum_{j=1}^n k_{i,j}e_j^2-\sigma_i^2\right\vert=o_{\P}\left(n^{-1/4}\right)$. 

\item
We have 
$$\max_{1\leq i\leq n}\left\vert\frac{1}{\hat{\sigma}_i^2}-\frac{1}{\sigma_i^2}\right\vert=o_{\P}\left(n^{-1/4}\right),\quad  \max_{1\leq i\leq n}\left\vert \frac{1}{\hat{\sigma}_i^2}-\frac{1}{\sigma_i^2}-\frac{\sigma_i^2-\sum_{j=1}^n k_{i,j}e_j^2}{\sigma_i^4}\right\vert=o_{\P}\left(n^{-1/2}\right).$$

\item
Define 
$$\check{s}_{1,i}=\sum_{j=1}^n k_{i,j}\frac{x_{1,j}y_j}{\sigma_j^2},\quad\overline{s}_{1,i}=\sum_{j=1}^n k_{i,j}\left[\sigma_j^2-\sum_{t=1}^n k_{j,t}e_t^2\right]\frac{x_{2,j}y_j}{\sigma_j^4},$$
with similar definition for $\check{s}_{2,i},\check{s}_{3,i}$ and $\overline{s}_{2,i},\overline{s}_{3,i}$.
Then we have for $\ell=1,2,3$,
$$\max_{1\leq i\leq n}\left\vert \hat{s}^{*}_{\ell,i}-\check{s}_{\ell,i}-\overline{s}_{\ell,i}\right\vert=o_{\P}\left(n^{-1/2}\right).$$
Moreover $\max_{1\leq i\leq n}\left\vert \overline{s}_{\ell,i}\right\vert=o_{\P}\left(n^{-1/4}\right)$ and $\max_{1\leq i\leq n}\left\vert \left(\hat{s}^{*}_{3,i}\right)^{-1}\right\vert=O_{\P}(1)$. 
\item
We set $m^{*}_i=\E\left(\check{s}_{3,i}\right)$ and for $\ell=1,2$,  $\widetilde{q}^{*}_{\ell,i}=\left(m^{*}_i\right)^{-1}\E\left(\check{s}_{\ell,i}\right)$. Then,
$$\max_{1\leq i\leq n}\left\vert \hat{q}^{*}_{\ell,i}-\widetilde{q}^{*}_{\ell,i}-\left(m^{*}_i\right)^{-1}\left[\E\check{s}_{3,i}-\hat{s}^{*}_{3,i}\right]\left(m^{*}_i\right)^{-1}
-\left(m^{*}_i\right)^{-1}\left[\hat{s}^{*}_{\ell,i}-\E\check{s}_{\ell,i}\right]\right\vert=o_{\P}\left(n^{-1/2}\right).$$
Moreover, if $s_{\ell,i}^{*}$ and $q_{\ell,i}^{*}$ are defined as $s_{\ell,i}$ and $q_{\ell,i}$ but with a division by $\sigma_i^2$ inside all the mathematical expectations,
we have also 
$$\max_{1\leq i\leq n}\left[\left\vert \hat{q}^{*}_{\ell,i}-q^{*}_{\ell,i}\right\vert+\left\vert \widetilde{q}^{*}_{\ell,i}-q^{*}_{\ell,i}\right\vert\right]=o_{\P}\left(n^{-1/4}\right).$$ 

\end{enumerate}
\end{lem}

\paragraph{Proof of Lemma \ref{plugin}}
\begin{enumerate}
\item
For the first point, we set $d_i=\sum_{j=1}^n k_{i,j}x_jx_j'$ and we use the decomposition  
$$\hat{\beta}(i/n)-\beta(i/n)=d_i^{-1}\sum_{j=1}^n k_{i,j}x_je_j+d_i^{-1}\sum_{j=1}^n k_{i,j}x_jx_j'\left[\beta(j/n)-\beta(i/n)\right].$$
As in the proof of Lemma \ref{approxquotient}, point $1$, we have $\max_{1\leq i\leq n}\left\vert d_i^{-1}\right\vert=O_{\P}(1)$. Then the result follows from the Lipschitz properties of $\beta$, Lemma \ref{uniforme}, the second point of the present lemma and our bandwidth conditions.
\item
For the second point, we use the fact that $\sum_{n\geq 1,1\leq i\leq n}\E\vert x_i\vert^2<\infty$ and from Lemma $A.1$ in \citet{ZW},
$$\max_{1\leq i\leq n}\left\vert \sum_{j=1}^n k_{i,j}\left[\vert x_j\vert^2-\E\vert x_j\vert^2\right]\right\vert=O_{\P}\left(\frac{1}{\sqrt{n}b}\right).$$

\item
We have 
\begin{eqnarray*}
&&\hat{\sigma}^2(i/n)-\sum_{j=1}^n k_{i,j}e_j^2\\
&=& \sum_{j=1}^n k_{i,j}\left\vert x_j'\left(\beta(j/n)-\hat{\beta}(i/n)\right)\right\vert^2-2\sum_{j=1}^n k_{i,j}e_j x_j'\left(\hat{\beta}(i/n)-\beta(j/n)\right)\\
&=& A_i-2B_i.
\end{eqnarray*}
We first deal with $A_i$. We have 
$$\left\vert A_i\right\vert\leq 2\sum_{j=1}^n k_{i,j}\vert x_j\vert^2\left[\left\vert \beta(j/n)-\beta(i/n)\right\vert^2+
\left\vert \hat{\beta}(i/n)-\beta(i/n)\right\vert^2\right].$$
From the two first points and using our bandwidth conditions, we have 
$$\max_{1\leq i\leq n}\left\vert A_i\right\vert=O_{\P}\left(b^2\right)+o_{\P}\left(n^{-1/2}\right)=o_{\P}\left(n^{-1/2}\right).$$
Next, we have $B_i=B_{i1}+B_{i2}$ with 
$$B_{i1}=\sum_{j=1}^nk_{i,j}e_jx_j'\left[\hat{\beta}(i/n)-\beta(i/n)\right],\quad B_{i2}=\sum_{j=1}^nk_{i,j}e_jx_j'\left[\beta(i/n)-\beta(j/n)\right].$$
Under our bandwidth condition $c<\frac{3}{4}-\frac{1}{q_2}$ and the assumption on $\Theta_{n,q_2}(L)$, we get from Lemma \ref{uniforme}   $\max_{1\leq i\leq n}\left\vert \sum_{j=1}^n k_{i,j}x_je_j\right\vert=o_{\P}\left(n^{-1/4}\right)$.  From point $1.$, we deduce that $\max_{1\leq i\leq n}\left\vert B_{i1}\right\vert=o_{\P}\left(n^{-1/2}\right)$.
Finally, we have $B_{i2}=b\sum_{j=1}^ne_j x_j'\widetilde{k}_{i,j}$ with $\widetilde{k}_{i,j}=k_{i,j}b^{-1}\left(\beta(i/n)-\beta(j/n)\right)$. One can note that the weights $\widetilde{k}_{i,j}$ satisfies similar properties as $k_{i,j}$. Indeed, we have $\max_{1\leq i,j\leq n}\left\vert\widetilde{k}_{i,j}\right\vert=O\left(\frac{1}{nb}\right)$ and 
$$\max_{1\leq i\leq n}\max_{1\leq j\leq n-1}\left\vert \widetilde{k}_{i,j}-\widetilde{k}_{i,j+1}\right\vert=O\left(\frac{1}{n^2b^2}\right).$$
Then, we claim that Lemma \ref{uniforme} can be applied by replacing the weights $k_{i,j}$ by the weights $\widetilde{k}_{i,j}$ because the two previous properties on the weights are sufficient to get this result (see also the proof of Lemma $A.3$ in \citet{ZW}). 
Using our bandwidth conditions, we easily get $\max_{1\leq i\leq n}\left\vert B_{i2}\right\vert=o_{\P}\left(n^{-1/2}\right)$.
This completes the proof of the first assertion. The second assertions easily follows Lemma \ref{uniforme}, the Lipschitz property of $\sigma$ and our bandwidth conditions.
\item
Using the previous point, we have 
$$\max_{1\leq i\leq n}\left\vert\frac{1}{\hat{\sigma}^2_i}-\frac{1}{\sigma_i^2}-\frac{\sigma_i^2-\hat{\sigma}^2_i}{\sigma_i^4}-\frac{\left(\hat{\sigma}_i^2-\sigma_i^2\right)^2}{\sigma_i^4\hat{\sigma}^2_i}\right\vert=o_{\P}\left(n^{-1/2}\right).$$
Using point $3$ and the positivity of $\sigma$, we get the result. 
\item
Setting $z_{1,i}=x_{2,i}y_i$, $z_{2,i}=x_{2,i}x_{1,i}'$ and $z_{3,i}=x_{2,i}x_{2,i}'$, we have
$$\hat{s}^{*}_{\ell,i}=\check{s}_{\ell,i}+\sum_{j=1}^n k_{i,j}z_{\ell,i}\left[\frac{1}{\hat{\sigma}^2_i}-\frac{1}{\sigma_i^2}\right].$$
Then the first assertion follows from the previous point of the present lemma. For the second assertion, we note that 
$$\max_{1\leq i\leq n}\left\vert \overline{s}_{\ell,i}\right\vert\leq \max_{1\leq i\leq n}\left\vert\sigma_i^2-\sum_{j=1}^n k_{i,j}e_j^2\right\vert$$
and the result follows from the point $3$. Finally, the last assertion follows from the first one, by using 
$$\max_{1\leq i\leq n}\left\vert\check{s}_{3,i}-\E\left(\sigma_i^{-2}x_{2,i}x_{2,i}'\right)\right\vert=o_{\P}(1),\quad \max_{1\leq i\leq n}\left\vert \E\left(\sigma_i^{-2}x_{2,i}x_{2,i}'\right)-\E\left(\sigma_i^{-2}x_2(i/n)x_2(i/n)'\right)\right\vert=O(1/n)$$
and assumption {\bf A4}.
\item
The first assertion follows from the previous point, using the same decomposition as in Lemma \ref{approxquotient} for the difference $\hat{q}^{*}_{\ell,i}-\widetilde{q}^{*}_{\ell,i}$. 
The other assertions can be shown as in Lemma \ref{approxquotient}, points $3$ and $4$.$\square$
\end{enumerate}

\subsection{Proof of Theorem \ref{central2}}
As in the proof of Theorem \ref{central}, we set 
$D_n=\sum_{i=1}^n \frac{1}{\hat{\sigma}^2_i}\hat{x}^{*}_{1,i}\left(\hat{x}^{*}_{1,i}\right)'$ and
we use the decomposition 
$$\hat{y}^{*}_i=\left(\hat{x}^{*}_{1,i}\right)'\beta_1+x_{2,i}'\hat{r}^{*}_i+e_i,\quad \hat{r}^{*}_i=\left(\hat{q}^{*}_{2,i}-q^{*}_{2,i}\right)\beta_1-\left(\hat{q}^{*}_{1,i}-q^{*}_{1,i}\right).$$
We get 
$$\hat{\beta}^{*}_1-\beta_1=\left(D_n^{*}\right)^{-1}\sum_{i=1}^n \frac{1}{\hat{\sigma}_i^{2}}\hat{x}^{*}_{1,i}\left[x_{2,i}'\hat{r}^{*}_i+e_i\right].$$
Using now the decomposition $\hat{x}^{*}_{1,i}=\widetilde{x}^{*}_{1,i}-\left(\hat{q}^{*}_{2,i}-q^{*}_{2,i}\right)'x_{2,i}$, with 
$\widetilde{x}^{*}_{1,i}=x_{1,i}-\left(q^{*}_{2,i}\right)'x_{2,i}$, we have
$$\hat{\beta}^{*}_1-\beta_1=\left(D_n^{*}\right)^{-1}\left[N^{*}_1-N^{*}_2+N^{*}_3-N^{*}_4\right],$$
where 
$$N^{*}_1=\sum_{i=1}^n \frac{1}{\hat{\sigma}_i^2}\widetilde{x}^{*}_{1,i} x_{2,i}'\hat{r}^{*}_i,\quad N^{*}_2=\sum_{i=1}^n \frac{1}{\hat{\sigma}_i^2}\left(\hat{q}^{*}_{2,i}-q^{*}_{2,i}\right)'x_{2,i}x_{2,i}'\hat{r}^{*}_i,$$
$$N^{*}_3=\sum_{i=1}^n \frac{1}{\hat{\sigma}_i^2}\widetilde{x}^{*}_{1,i}e_i,\quad N^{*}_4=\sum_{i=1}^n\frac{1}{\hat{\sigma}_i^2}\left(\hat{q}^{*}_{2,i}-q^{*}_{2,i}\right)'x_{2,i}e_i.$$
\begin{itemize}
\item
We first consider the leading term $\frac{1}{\sqrt{n}}N^{*}_3$. Using Lemma \ref{plugin}, we have
$$\frac{1}{\sqrt{n}}N^{*}_3=\frac{1}{\sqrt{n}}\sum_{i=1}^n \frac{\widetilde{x}_{1,i}^{*}e_i}{\sigma_i^2}+\frac{1}{\sqrt{n}}\sum_{i=1}^n \frac{\sigma_i^2-\sum_{j=1}^n k_{i,j}e_j^2}{\sigma_i^4}\widetilde{x}^{*}_{1,i}e_i+o_{\P}(1).$$ 
As in the proof of Theorem \ref{central}, one can show that the first term is asymptotically Gaussian with mean $0$ and variance $\Sigma^{*}$. Note that one can use directly the central limit theorem for martingale differences for this step.
For the second term, we decompose 
$$\sum_{i=1}^n \frac{\sigma_i^2-\sum_{j=1}^n k_{i,j}e_j^2}{\sigma_i^4}\widetilde{x}^{*}_{1,i}e_i=S_1+S_2,$$
with $S_1=\sum_{i,j=1}^n k_{i,j}\left(\sigma_j^2-\sigma_i^2\right)\frac{\widetilde{x}^{*}_{1,i}e_i}{\sigma_i^4}=o_{P}\left(\sqrt{n}\right)$ by using the moment inequality given in \citet{ZW}, Lemma $A.1$ and 
$S_2=\sum_{i,j=1}^n k_{i,j}\left(e_j^2-\sigma_j^2\right)\frac{\widetilde{x}^{*}_{1,i}e_i}{\sigma_i^4}=o_{P}\left(\sqrt{n}\right)$
by using Lemma \ref{smooth1} with $r=4/3$, $r_1=2$ and $r_2=4$.  
\item
From Lemma \ref{plugin}, we easily get $N^{*}_2=o_{\P}\left(\sqrt{n}\right)$. 
\item
Next, we show that $N^{*}_4=o_{\P}\left(\sqrt{n}\right)$. Using Lemma \ref{plugin}, it is easily seen that 
$$N^{*}_4=\sum_{i=1}^n\left(\hat{q}^{*}_{2,i}-q^{*}_{2,i}\right)'\frac{x_{2,i}e_i}{\sigma_i^2}+o_{\P}\left(\sqrt{n}\right).$$
Moreover, using Lemma $A.1$ in \citet{ZW}, we have 
$$N^{*}_4=\sum_{i=1}^n \left(\hat{q}^{*}_{2,i}-\widetilde{q}^{*}_{2,i}\right)'\frac{x_{2,i}e_i}{\sigma_i^2}+o_{\P}\left(\sqrt{n}\right).$$
Using Lemma \ref{plugin}$(6)$, it remains to show that 
$\sum_{i=1}^n z_i'\frac{x_{2,i}e_i}{\sigma_i^2}=o_{\P}\left(\sqrt{n}\right)$ when 
$$z_i=\left(m^{*}_i\right)^{-1}\left[m_i-\hat{s}^{*}_{3,i}\right]\left(m^{*}_i\right)^{-1}\mbox{ or } z_i=\left(m^{*}_i\right)^{-1}\left[\hat{s}^{*}_{2,i}-\E\check{s}_{2,i}\right].$$
Let us detail the first case. We have $m_i^{*}z_im_i^{*}=z_{1,i}+z_{2,i}$ with $z_{1,i}=\check{s}_{3,i}-\E\check{s}_{3,i}$ and $z_{2,i}=\overline{s}_{3,i}$. For $z_{1,i}$, one can use Lemma \ref{smooth1} with $r_1=q_1$ and $r_2=q_2$. For the $z_{2,i}$, we have 
$$\left\vert \sum_{i=1}^n \left(m^{*}_i\right)^{-1}z_{2,i}'\left(m^{*}_i\right)^{-1}\frac{x_{2,i}e_i}{\sigma_i^2} \right\vert\leq Cn\max_{1\leq j\leq n}\left\vert \sum_{i=1}^n k_{i,j}\frac{x_{2,i}'e_i}{\sigma_i^2}\right\vert\cdot\max_{1\leq j\leq n}\left\vert \sum_{t=1}^n k_{j,t}e_t^2-\sigma_j^2\right\vert.$$
Then the result follows from Lemma \ref{uniforme}. Note that the first maximum involves the weights $k_{i,j}$ instead of $k_{j,i}$ but it does not change the result and this maximum is a $o_{\P}\left(n^{-1/4}\right)$. Moreover, the second maximum is $o_{\P}\left(n^{-1/4}\right)$ by applying Lemma \ref{plugin}$(3)$.
The case $z_i=\left(m^{*}_i\right)^{-1}\left[\hat{s}^{*}_{2,i}-\E\check{s}_{2,i}\right]$ is similar and then omitted. 
Using analog arguments, one can also show that $N^{*}_1=o_{\P}\left(\sqrt{n}\right)$. Details are omitted.

\item
Finally, we show that $D^{*}_n/n\rightarrow \Sigma^{*}$ in probability. Using the decomposition
\begin{eqnarray*}
D^{*}_n&=&\sum_{i=1}^n \frac{1}{\hat{\sigma}^2_i}\widetilde{x}^{*}_{1,i}\left(\widetilde{x}^{*}_{1,i}\right)'+\sum_{i=1}^n\left(\hat{q}^{*}_{2,i}-q^{*}_{2,i}\right)\frac{x_{2,i}x_{2,i}'}{\hat{\sigma}^2_i}\left(\hat{q}^{*}_{2,i}-q^{*}_{2,i}\right)'\\
&+& \sum_{i=1}^n\left(\hat{q}^{*}_{2,i}-q^{*}_{2,i}\right)' \frac{x_{2,i}\widetilde{x}_{1,i}'}{\hat{\sigma}_i^2}+\sum_{i=1}^n \frac{\widetilde{x}_{1,i}x_{2,i}'}{\hat{\sigma}_i^2}\left(\hat{q}^{*}_{2,i}-q^{*}_{2,i}\right).
\end{eqnarray*}
Using the previous points and Lemma \ref{plugin}$(4)$, we get 
$$\frac{D^{*}_n}{n}=\frac{1}{n}\sum_{i=1}^n \frac{\widetilde{x}^{*}_{1,i}\left(\widetilde{x}^{*}_{1,i}\right)'}{\hat{\sigma}_i^2}+o_{\P}(1)=\frac{1}{n}\sum_{i=1}^n \frac{\widetilde{x}^{*}_{1,i}\left(\widetilde{x}^{*}_{1,i}\right)'}{\sigma_i^2}=o_{\P}(1).$$
As for the end of the proof of Theorem \ref{central}, we get $n^{-1}D^{*}_n=\Sigma^{*}+o_{\P}(1).$
$\square$

\end{itemize}
\subsection{Auxiliary results}

We first recall a lemma given in \citet{ZW} for uniform convergence of kernel smoothers.

\begin{lem}\label{uniforme}
Let $z_{n,i}=J_{n,i}\left(\mathcal{F}_i\right)$ be a centered triangular array. We assume that there exists $r>2$ such that $\sup_{1\leq i\leq n,n\geq 1}\Vert z_{n,i}\Vert_r<\infty$ and $\Theta_{n,r}(J)=O\left(n^{-\nu}\right)$ for $\nu>\frac{1}{2}-\frac{1}{r}$. Then,  
$$\sup_{u\in [0,1]}\left\vert \sum_{i=1}^n \frac{1}{nb}K\left(\frac{u-i/n}{b}\right)z_{n,i}\right\vert=\frac{O_{\P}(v_n)}{nb},\mbox{ where } v_n=n^{1/r}+\sqrt{nb \log n}.$$
\end{lem} 
The second lemma is useful for getting convergence rates of quadratic forms.

\begin{lem}\label{controlquad}
Let $z_{n,i}=J_{1,n,i}\left(\mathcal{F}_i\right)$ and $w_{n,i}=J_{2,n,i}\left(\mathcal{F}_i\right)$ be two centered triangular arrays.
Let $q\in (1,2]$ and $q_1,q_2>2$ two real numbers such that $\frac{1}{q_1}+\frac{1}{q_2}=\frac{1}{q}$ and
$\Theta_{0,q_1}\left(J_1\right)+\Theta_{0,q_2}\left(J_2\right)<\infty$ and $Q_n=\sum_{1\leq i<j\leq n}\alpha_{i,j}z_{n,i}w_{n,j}$, where for $1\leq i<j\leq n$, $\alpha_{i,j}$ is a real number. Then we have 
$$\Vert Q_n-\E Q_n\Vert_q\leq C n^{1/q}\left[\Theta_{0,q_1}(J_1)+\Theta_{0,q_2}(J_2)\right]^2a_n,$$
where
$$a_n^2=\left(\max_{1\leq i<n}\sum_{j=i+1}^n\vert\alpha_{i,j}\vert^2\right)\vee\left(\max_{1<j\leq n}\sum_{i=1}^{j-1}\vert \alpha_{i,j}\vert^2\right).$$ 
\end{lem}

\paragraph{Proof of Lemma \ref{controlquad}}
We follows the proof of Proposition $1$ in \citet{LW} which gives a control of $\L^q$ norms of quadratic forms when $q\geq 2$. 
See also Lemma $A.2$ in \citet{ZW} for a similar result.
We set 
$$u_{n,j}=\sum_{i=1}^{j-1}\alpha_{i,j} z_{n,i},\quad v^{(k)}_{n,i}=\sum_{j=i+1}^n \alpha_{i,j}w_{n,j}^{(k)},$$
where for $k\in\Z$, $w_{n,j}^{(k)}$ is defined as $w_{n,j}$ but the innovation $\xi_k$ is replaced with $\xi_k'$ in $\mathcal{F}_j$ ($\xi'$ is a copy of $\xi$).
Defining $$\mathcal{P}_k=\E\left(\cdot\vert \mathcal{F}_k\right)-\E\left(\cdot\vert \mathcal{F}_{k-1}\right),$$
we have 
\begin{eqnarray*}
\Vert \mathcal{P}_k Q_n\Vert_q &\leq& \Vert \sum_{1\leq i<j\leq n} \alpha_{i,j} \left[z_{n,i}w_{n,j}-z_{n,i}^{(k)}w_{n,j}^{(k)}\right]\Vert_q\\
&\leq& \Vert \sum_{j=2}^n u_{n,j}\left[w_{n,j}-w_{n,j}^{(k)}\right]\Vert_q+\Vert \sum_{i=1}^{n-1}v^{(k)}_{n,i}\left[z_{n,i}-z_{n,i}^{(k)}\right]\Vert_q\\
&=& I_k+II_k.
\end{eqnarray*} 
Using the moment bound given in Lemma $A.1$ in \citet{ZW}, we have 
\begin{eqnarray*}
I_k&\leq& \sum_{j=2}^n \Vert u_{n,j}\Vert_{q_1}\cdot \Vert w_{n,j}-w_{n,j}^{(k)}\Vert_{q_2}\\
&\leq& Ca_n\sum_{j=2}^n \delta_{j-k,q_2}\left(J_2\right)\Theta_{0,q_1}\left(J_1\right). 
\end{eqnarray*}
Using similar arguments, we get 
$$II_k\leq C a_n \sum_{i=1}^{n-1}\delta_{i-k,q_1}\left(J_1\right) \Theta_{0,q_2}\left(J_2\right).$$
Since 
$$\sum_{j=2}^n\delta_{j-k,2q}(J_2)\leq \Theta_{0,q_2}(J_2),\quad \sum_{i=1}^{n-1}\delta_{i-k,2q}(J_1)\leq \Theta_{0,q_1}(J_1),$$
we get 
\begin{eqnarray*}
\Vert \mathcal{P}_k Q_n\Vert_q^q &\leq & C\left[I_k^q+II_k^q\right]\\
&\leq& C a_n^q \left(\Theta_{0,q_1}(J_1)+\Theta_{0,q_2}(J_2)\right)^{2q-1}\left[\sum_{j=2}^n \delta_{j-k,2q}(J_2)+\sum_{i=1}^{n-1}\delta_{i-k,2q}(J_1)\right].
\end{eqnarray*}
Finally, using Burkholder's inequality, we get 
\begin{eqnarray*}
\Vert Q_n-\E Q_n\Vert_q&=&\Vert \sum_{k=-\infty}^n \mathcal{P}_k Q_n\Vert_q\\
&\leq& C \E^{1/q}\left[\left\vert \sum_{k=-\infty}^n \left[\mathcal{P}_k Q_n\right]^2\right\vert^{q/2}\right]\\
&\leq& C \left[\sum_{k=-\infty}^n \E\left\vert \mathcal{P}_k Q_n\right\vert^q\right]^{1/q}\\
&\leq & C a_n n^{1/q}\left[\Theta_{0,q_1}(J_1)+\Theta_{0,q_2}(J_2)\right]^2.
\end{eqnarray*}
This completes the proof of the lemma.$\square$

\begin{lem}\label{smooth1}
Let $z_i=J_{1,n,i}\left(\mathcal{F}_i\right)$ and $w_i=J_{2,n,i}\left(\mathcal{F}_i\right)$ be two centered triangular arrays. Assume that there exists $r\in (1,2]$ such that $\Theta_{0,r_1}(J_1)+\Theta_{0,r_2}(J_2)<\infty$ where $r_1,r_2>2$ are such that $\frac{1}{r}=\frac{1}{r_1}+\frac{1}{r_2}$. Then, if $n^db\rightarrow \infty$ for $d=\min\left(2-\frac{2}{r},\frac{1}{2}\right)$, 
$$\frac{1}{\sqrt{n}}\sum_{i=1}^nz_i\cdot\sum_{j=1}^nk_{i,j}w_i=o_{\P}(1).$$
\end{lem}

\paragraph{Proof of Lemma \ref{smooth1}}
Since $\max_{1\leq i\leq n}\E\left\vert z_i w_i\right\vert=O(1)$, we have 
$$\frac{1}{\sqrt{n}}\sum_{i=1}^n k_{i,i}z_i w_i=O_{\P}\left(\frac{1}{\sqrt{n}b}\right)=o_{\P}(1).$$
Moreover, we have 
$$\max_{1\leq i<n}\left(\sum_{j=i+1}^nk^2_{i,j}\right)\vee\left(\max_{1<j\leq n}\sum_{i=1}^{j-1} k_{i,j}^2\right)=O\left(\frac{1}{nb}\right).$$
Using Lemma \ref{controlquad}, we get under our bandwidth conditions 
$$\frac{1}{\sqrt{n}}\sum_{1\leq i\neq j\leq n} k_{i,j}z_i w_j=O_{\P}\left(\frac{n^{1/r}}{n\sqrt{b}}\right)=o_{\P}(1).$$
Finally, we have 
\begin{equation}\label{finalstep}
\frac{1}{\sqrt{n}}\sum_{1\leq i\neq j\leq n}k_{i,j}\E\left(z_iw_j\right)=O\left(\frac{1}{\sqrt{n}b}\right). 
\end{equation}
To show (\ref{finalstep}), we bound the covariances as follows. If $i\leq j$, we have 
\begin{eqnarray*}
\left\vert \E\left(z_i w_j\right)\right\vert&\leq& \sum_{\ell\leq i}\left\vert \E\left(\mathcal{P}_{\ell}z_i\cdot\mathcal{P}_{\ell}w_j\right)\right\vert\\
&\leq& \sum_{\ell\leq i}\delta_{i-\ell,r_1}(J_1)\delta_{j-\ell,r_2}(J_2).
\end{eqnarray*}
Hence, we obtain 
$$\frac{1}{\sqrt{n}}\sum_{1\leq i,j\leq n}k_{i,j}\left\vert\E\left(z_iw_j\right)\right\vert\leq \sqrt{n}\max_{1\leq i,j\leq n}k_{i,j}\Theta_{0,r_1}(J_1)\Theta_{0,r_2}(J_2),$$
which leads to (\ref{finalstep}) since $\max_{1\leq i,j\leq n}k_{i,j}=O\left(\frac{1}{nb}\right)$. The proof of Lemma \ref{smooth1}
is now complete.$\square$

\begin{lem}\label{util1}
Let $\Sigma:[0,1]\rightarrow \mathcal{S}^{+}_d$ be a continuous application taking values in the space $\mathcal{S}^{+}_d$ of square matrices of size $d$, symmetric and positive definite. Assume that for each $u\in[0,1]$, $\Sigma(u)$ is invertible. Then we have the inequality 
$$\left(\int_0^1\Sigma(u)du\right)^{-1}\preceq \int_0^1\Sigma(u)^{-1}du,$$
where $\preceq$ denotes the order relation for nonnegative definite matrices. Moreover, the equality between these two matrices holds if and only if the application $\Sigma$ is constant.  
\end{lem}

\paragraph{Proof of Lemma \ref{util1}}

For $x,y\in \R^d$, we have from the Cauchy-Schwarz inequality, 
$$(x'\cdot y)^2=\left\vert \int_0^1 x'\cdot\Sigma(u)^{-1/2}\Sigma(u)^{1/2}\cdot y du\right\vert^2\leq x'\cdot \int_0^1 \Sigma(u)^{-1}du \cdot x\cdot y'\cdot\int_0^1\Sigma(u)du\cdot y.$$
Then setting $y=\left(\int_0^1 \Sigma(u)du\right)^{-1}\cdot x$, we get 
$$x'\cdot\left(\int_0^1\Sigma(u)du\right)^{-1}\cdot x\leq x'\cdot\int_0^1\Sigma(u)^{-1}du\cdot x,$$
which shows the inequality. Moreover, if the two matrices are equal, the equality condition in the Cauchy-Schwarz inequality entails that for each $x\in\R^d\setminus \{0\}$, there exists a real number $\lambda_x$ such that $x=\lambda_x \Sigma(u)\left(\int_0^1\Sigma(u)du\right)^{-1}x$. Integrating the previous relation with respect to $u$, we find that $\lambda_x=1$ and 
then $\Sigma(v)=\int_0^1\Sigma(u)du$, $v\in [0,1]$.$\square$

\bibliographystyle{plainnat}
\bibliography{bibpartial}
\end{document}